\documentclass[openany, a4paper, oneside, 12pt]{myamsbook}
\usepackage{amsmath,amssymb, amsfonts}
\usepackage[latin10]{inputenc}
\usepackage[english]{babel}
\usepackage{amscd}
\usepackage{epsfig, epsf, color, xypic}
\usepackage[all,cmtip]{xy}

\usepackage{fancyhdr}

\usepackage{fourier}

\usepackage[dvips]{geometry}

\renewcommand{\chaptermark}[1]{\markboth{#1}{}}

\begin{document}

\pagenumbering{arabic}

\newcommand{\captionfonts}{\small}

\makeatletter  %
\long\def\@makecaption#1#2{%
  \vskip\abovecaptionskip
  \sbox\@tempboxa{{\captionfonts #1: #2}}%
  \ifdim \wd\@tempboxa >\hsize
    {\captionfonts #1: #2\par}
  \else
    \hbox to\hsize{\hfil\box\@tempboxa\hfil}%
  \fi
  \vskip\belowcaptionskip}

\newcommand*\if@single[3]{%
  \setbox0\hbox{${\mathaccent"0362{#1}}^H$}%
  \setbox2\hbox{${\mathaccent"0362{\kern0pt#1}}^H$}%
  \ifdim\ht0=\ht2 #3\else #2\fi
  }
\newcommand*\rel@kern[1]{\kern#1\dimexpr\macc@kerna}
\newcommand*\widebar[1]{\@ifnextchar^{{\wide@bar{#1}{0}}}{\wide@bar{#1}{1}}}
\newcommand*\wide@bar[2]{\if@single{#1}{\wide@bar@{#1}{#2}{1}}{\wide@bar@{#1}{#2}{2}}}
\newcommand*\wide@bar@[3]{%
  \begingroup
  \def\mathaccent##1##2{%
    \if#32 \let\macc@nucleus\first@char \fi
    \setbox\z@\hbox{$\macc@style{\macc@nucleus}_{}$}%
    \setbox\tw@\hbox{$\macc@style{\macc@nucleus}{}_{}$}%
    \dimen@\wd\tw@
    \advance\dimen@-\wd\z@
    \divide\dimen@ 3
    \@tempdima\wd\tw@
    \advance\@tempdima-\scriptspace
    \divide\@tempdima 10
    \advance\dimen@-\@tempdima
    \ifdim\dimen@>\z@ \dimen@0pt\fi
    \rel@kern{0.6}\kern-\dimen@
    \if#31
      \overline{\rel@kern{-0.6}\kern\dimen@\macc@nucleus\rel@kern{0.4}\kern\dimen@}%
      \advance\dimen@0.4\dimexpr\macc@kerna
      \let\final@kern#2%
      \ifdim\dimen@<\z@ \let\final@kern1\fi
      \if\final@kern1 \kern-\dimen@\fi
    \else
      \overline{\rel@kern{-0.6}\kern\dimen@#1}%
    \fi
  }%
  \macc@depth\@ne
  \let\math@bgroup\@empty \let\math@egroup\macc@set@skewchar
  \mathsurround\z@ \frozen@everymath{\mathgroup\macc@group\relax}%
  \macc@set@skewchar\relax
  \let\mathaccentV\macc@nested@a
  \if#31
    \macc@nested@a\relax111{#1}%
  \else
    \def\gobble@till@marker##1\endmarker{}%
    \futurelet\first@char\gobble@till@marker#1\endmarker
    \ifcat\noexpand\first@char A\else
      \def\first@char{}%
    \fi
    \macc@nested@a\relax111{\first@char}%
  \fi
  \endgroup
}

\makeatother   %

\newcommand{\ii}{\^{i}}
\newcommand{\Ii}{\^{I}}
\newcommand{\sz}{\c{s}}
\newcommand{\ty}{\c{t}}
\newcommand{\ti}{\ty i}
\newcommand{\aw}{\u{a}}
\newcommand{\io}{\^{a}}

\newcommand{\naX}{\nabla_X}
\newcommand{\na}{\nabla}
\newcommand{\Ti}{T_{e_i}}
\newcommand{\Tj}{T_{e_j}}

\newcommand{\jjg}{,\hspace{-1.5pt},}
\newcommand{\si}{\c{s}i}
\newcommand{\cn}{c\io nd}
\newcommand{\nn}{\mbox {\boldmath $n$}}
\newcommand{\tg}{\mbox {\boldmath $t$}}
\newcommand{\bb}{\mbox {\boldmath $b$}}
\newcommand{\rg}{\mbox {\rm rg}\,}
\newcommand{\grad}{\mathop\mathrm{grad} }
\newcommand{\fa}{func\c{t}ia}
\newcommand{\fe}{func\ty ie}
\newcommand{\Fa}{Func\c{t}ia}
\newcommand{\su}{suprafa\c{t}\u{a}}
\newcommand{\sa}{S\u{a} se ar\aw te c\u{a}}
\newcommand{\ve}{vecin\u{a}tate}
\newcommand{\va}{variet\u{a}\c{t}i}
\newcommand{\so}{\bf Solu\c{t}ie}
\newcommand{\dd}{dac\u{a} \sz i\ numai dac\u{a}}
\newcommand{\ec}{ecua\c{t}i}
\newcommand{\ai}{astfel \^{\i}nc\io t}
\newcommand{\ci}{c\^{\i}mp}
\newcommand{\diff}{diferen\ty abil}
\newcommand{\s}{\sigma}

\newcommand{\Lg}{\mathfrak{g}}

\newcommand{\K}{\mathop\mathrm{Ker}}
\newcommand{\I}{\mathop\mathrm{Im}}
\newcommand{\norm}[1]{\Vert #1\Vert}
\newcommand{\inn}[1]{#1_\infty}
\newcommand{\abs}[1]{\vert #1\vert}
\newcommand{\fr}[2]{\frac{\vv{#1}}{\vv{#2}}}
\newcommand{\ov}{\overline}
\newcommand{\vv}{\overrightarrow}
\newcommand{\6}[1]{\frac{\partial}{\partial x^{#1}}}
\newcommand{\tl}[1]{\tilde{#1}}
\newcommand{\Tl}[1]{\widetilde{#1}}
\newcommand{\fra}{\frac{\mathop\mathrm{d}}{\mathop\mathrm{dt}}|_0}
\newcommand{\pa}[2]{\frac{\partial #1}{\partial #2}}
\newcommand{\dpa}[1]{\frac{D}{\partial #1}}
\newcommand{\ser}[2]{\int_{#1}^{#2}}
\newcommand{\suma}[2]{\sum_{#1}^{#2}}
\newcommand{\cp}{\bar{\partial}}
\newcommand{\cz}{\widebar{z}}

\newcommand{\Int}{\mathrm{Int}}
\newcommand{\tr}{\mathrm{tr}}
\newcommand{\Hom}{\mathrm{Hom}}

\newcommand{\zz}{\mathfrak{z}}
\newcommand{\x}{\mathcal{X}}
\newcommand{\ic}{\mathcal{I}}
\newcommand{\ce}{\mathcal{C}^\infty}
\newcommand{\A}{\mathcal{A}}
\newcommand{\B}{\mathcal{B}}
\newcommand{\C}{\mathcal{C}}
\newcommand{\D}{\mathcal{D}}
\newcommand{\E}{\mathcal{E}}
\newcommand{\F}{\mathcal{F}}
\newcommand{\Ll}{\mathfrak{L}}
\newcommand{\M}{\mathcal{M}_{k,n}}
\newcommand{\Oo}{\mathcal{O}}
\newcommand{\R}{\mathcal{R}}
\newcommand{\Ss}{\mathcal{S}}
\newcommand{\TT}{\mathbb{T}}
\newcommand{\HA}{\mathcal{H}}
\newcommand{\V}{\mathcal{V}}
\newcommand{\W}{\mathcal{W}}
\newcommand{\Ix}[1]{\mathcal{I}_{#1}}
\newcommand{\f}{\varphi}
\newcommand{\al}{\alpha}
\newcommand{\be}{\beta}
\newcommand{\g}{\gamma}
\newcommand{\G}{\Gamma}
\newcommand{\e}{\varepsilon}
\newcommand{\la}{\lambda}
\newcommand{\sgn}{\mathrm{sgn}}
\newcommand{\ra}{\rightarrow}
\newcommand{\gr}{\mathrm{gr}}
\newcommand{\GL}{\mathrm{GL}}
\newcommand{\Ric}{\textnormal{Ric}}

\newcommand{\ad}{\mathrm{ad}}
\newcommand{\Ad}{\mathrm{Ad}}

\newcommand{\dist}{\mathop\mathrm{d}}
\newcommand{\diam}{\mathop\mathrm{diam}}
\renewcommand{\tan}{\mathop\mathrm{tg}}
\renewcommand{\arctan}{\mathop\mathrm{arctg}}
\renewcommand{\cosh}{\mathop\mathrm{ch }}
\renewcommand{\sinh}{\mathop\mathrm{sh}}
\renewcommand{\cot}{\mathop\mathrm{ctg}}
\renewcommand{\div}{\mathop\mathrm{div }}
\newcommand{\rot}{\mathop\mathrm{rot}}
\newcommand{\vol}{\mathop\mathrm{vol}}
\renewcommand{\phi}{\varphi}

\newcommand{\CC}{\mathbb{C}}
\newcommand{\HH}{\mathbb{H}}
\newcommand{\RR}{\mathbb{R}}
\newcommand{\KK}{\mathbb{K}}
\newcommand{\ZZ}{\mathbb{Z}}
\newcommand{\NN}{\mathbb{N}}
\newcommand{\PP}{\mathbb{P}}
\newcommand{\Pn}[1]{\mathbb{P}^{#1}}
\newcommand{\QQ}{\mathbb{Q}}
\newcommand{\Gkn}[2]{G_{#1,#2}}
\newcommand{\Rn}{\RR^{2n}}
\newcommand{\Cn}[1]{\CC^{#1}}
\newcommand{\ip}{\frac{i}{2\pi}}
\newcommand{\overbar}[1]{\mkern 1.5mu\overline{\mkern-1.5mu#1\mkern-1.5mu}\mkern 1.5mu}
\newcommand{\Ka}{K\"{a}hler}
\newcommand{\hK}{hyperk\"{a}hler}
\newcommand{\HK}{Hyperk\"{a}hler}
\newcommand{\Tu}{T^{1,0}M}
\newcommand{\Td}{T^{0,1}M}
\newcommand{\PC}{\Pn{n}\CC}
\newcommand{\pl}{\square}
\newcommand{\cpl}{\widebar{\pl}}
\newcommand{\Lt}{\tl{L}}
\newcommand{\PN}{\Pn{N}\CC}

\numberwithin{equation}{section}
\numberwithin{figure}{chapter}

\newenvironment{demo}{\noindent{\bf Proof.}\hspace*{0.3em}}{{$\blacksquare$}\\}

\newenvironment{demon}{\noindent{\bf Demonstra\c
\ty ie Teorema 4.2.1.}\hspace*{2em}}{{$\blacksquare$}\\}

\newenvironment{demoK}{\noindent{\bf Demonstra\c
\ty ie la teorema de scufundare: }\hspace*{2em}}{{$\blacksquare$}\\}

\newenvironment{demo1}{\noindent{\bf Proof of Proposition 4.1.}\hspace*{2em}}{{$\blacksquare$}\\}

\renewcommand{\chaptername}{\bf Chapter}

\renewcommand{\contentsname}{Contents}
\renewcommand{\bibname}{Bibliography}
\renewcommand{\figurename}{Figura}
\renewcommand{\tablename}{Tabelul}
\renewcommand{\indexname}{\Ii ndice de no\c tiuni}

\addto\captionsromanian{\renewcommand{\partname}{}}

\newcommand{\partea}[1]{\let\oldthepart\thepart
  \renewcommand{\thepart}{\oldthepart}
  \part{#1}\let\thepart\oldthepart}

\newcommand{\parteaa}[1]{\let\oldthepart\thepart
  \renewcommand{\thepart}{a \oldthepart-a}
  \part{#1}\let\thepart\oldthepart}

\newcounter{Mycounter}[section]

\newcounter{lemma}[section]
\setcounter{lemma}{0}
\renewcommand{\thelemma}{{Lemma \thesection.\arabic{lemma}}}
\newcommand{\lm}{%
     \setcounter{lemma}{\value{Mycounter}}
     \refstepcounter{lemma}
     \stepcounter{Mycounter}
     {\vspace{1mm}\noindent \bf \thelemma.\ }\em }

\newcounter{corollary}[section]
\setcounter{corollary}{0}
\renewcommand{\thecorollary}{{Corollary
\thesection.\arabic{corollary}}}
\newcommand{\co}{%
     \setcounter{corollary}{\value{Mycounter}}
     \refstepcounter{corollary}
     \stepcounter{Mycounter}
     {\vspace{1mm}\noindent \bf \thecorollary.\ }\em }

\newcounter{theorem}[section]
\setcounter{theorem}{0}
\renewcommand{\thetheorem}{{Theorem \thesection.\arabic{theorem}}}
\newcommand{\te}{%
     \setcounter{theorem}{\value{Mycounter}}
     \refstepcounter{theorem}
     \stepcounter{Mycounter}
     {\vspace{1mm}\noindent \bf \thetheorem.\ }\em }

\newcounter{proposition}[section]
\setcounter{proposition}{0}
\renewcommand{\theproposition}{{Proposition \thesection.\arabic{proposition}}}
\newcommand{\pr}{%
     \setcounter{proposition}{\value{Mycounter}}
     \refstepcounter{proposition}
     \stepcounter{Mycounter}
     {\vspace{1mm}\noindent \bf \theproposition.\  }\em }

\newcounter{definition}[section]
\setcounter{definition}{0}
\renewcommand{\thedefinition}
       {{Definition~\thesection.\arabic{definition}}}
\newcommand{\de}{%
     \setcounter{definition}{\value{Mycounter}}
     \refstepcounter{definition}
     \stepcounter{Mycounter}
     {\vspace{1mm}\noindent \bf \thedefinition.\ }}

\newcounter{example}[section]
\setcounter{example}{0}
\renewcommand{\theexample}{{Example \thesection.\arabic{example}}}
\newcommand{\ex}{%
     \setcounter{example}{\value{Mycounter}}
     \refstepcounter{example}
     \stepcounter{Mycounter}
     {\vspace{1mm}\noindent \bf \theexample.\ }}
     
\newcounter{exercice}[section]
\setcounter{exercice}{0}
\renewcommand{\theexercice}{{Exercice \thesection.\arabic{exercice}}}
\newcommand{\apl}{%
     \setcounter{exercice}{\value{Mycounter}}
     \refstepcounter{exercice}
     \stepcounter{Mycounter}
     {\vspace{1mm}\noindent \bf \theexercice.\ }}

\newcounter{remark}[section]
\setcounter{remark}{0}
\renewcommand{\theremark}{{Remark \thesection.\arabic{remark}}}
\newcommand{\ob}{%
     \setcounter{remark}{\value{Mycounter}}
     \refstepcounter{remark}
     \stepcounter{Mycounter}
     {\vspace{1mm}\noindent \bf \theremark.\ }}

\newcounter{problem}[section]
\setcounter{problem}{0}
\renewcommand{\theproblem}{{{\small Exerci\c tiul}
\small\thechapter.\thesection.\arabic{problem}}}
\newcommand{\exc}{%
    \setcounter{problem}{\value{Mycounter}}
   \refstepcounter{problem}
   \stepcounter{Mycounter}
   {\noindent \bf \theproblem.\ }}

\mainmatter

\begin{center}
 \textbf{ TWISTED HOLOMORPHIC SYMPLECTIC FORMS}

\vspace {5pt}
Nicolina Istrati
\end{center}

\bigskip
\begin{quote}
{\tiny ABSTRACT. We show that a compact \Ka\ manifold admitting a nondegenerate holomorphic 2-form valued in a line bundle is a finite cyclic cover of a \hK\ manifold. With respect to the connection induced by the locally \hK\ metric, the form is parallel. We then describe the structure of the fundamenal group of such manifolds and derive some consequences.}

\end{quote}

\section{Introduction}

In this note, we are concerned with compact \Ka\ manifolds which admit a particular kind of structure: holomorphic nondegenerate 2-forms valued in a line bundle. The problem has different analogues that have been intensively studied. On the one hand, there is the nontwisted problem concerning holomorphic symplectic forms. On the other hand, its symmetric avatar consists in the study of holomorphic (conformal) metrics.

In the compact setting, \Ka\ manifolds admitting holomorphic symplectic forms are exactly the \hK\ ones, as shown in \cite{b}. There is a rich literature concerning this subject, and its study is ongoing. Turning to the symmetric counterpart, the situation is somewhat different. Although the class of compact \Ka\ manifolds admitting a holomorphic metric is rather small (it consists in all the finite coverings of tori, as shown in \cite{iko}), as soon as one allows the same structure to be twisted -- thus studying holomorphic conformal structures -- one enters a very rich class of manifolds. A complete classification of these has been reached only in dimension 2 and 3, in \cite{ko} and \cite{jr}.

Despite one could expect that the class of manifolds with twisted holomorphic symplectic forms is also wide, it turns out that the situation is not much different from the nontwisted case. More precisely, we show that compact \Ka\ manifolds admitting a twisted holomorphic symplectic form are locally \hK. In particular, the presence of such a structure ensures the existence of a Ricci-flat \Ka\ metric, and with respect to the connection induced by this metric the form is parallel.

Roughly speaking, the proof is as follows: we first notice that the twisted holomorphic symplectic form induces local Lefschetz-type operators acting on the sheaves of holomorphic forms $\Omega^*$, which then determine a local splitting of $\Omega^3$ into $\Omega^1$ and some other summand. This, in turn, allows us to find local holomorphic 1-forms which behave like connection forms on the line bundle where the twisted form takes its values. Finally, this means that the bundle admits a holomorphic connection, thus also a flat one, and that the manifold is Ricci-flat locally holomorphic symplectic, thus locally \hK. This is Theorem 2.7 in Section 2.

In the next section, we investigate the conditions that make the converse true, i.e. when does a locally \hK\ manifold admit a twisted holomorphic symplectic form. This is strictly related to the structure of its fundamental group: it has to be cyclic and of a particular kind. There are some consequences that follow, such as: (strictly) twisted holomorphic symplectic manifolds are necessarily irreducible, and if their fundamental group is finite, they are necessarily projective.

\bigskip

Let us now be more precise, and define the objects we will be interested in:

\begin{de}
A Riemannian manifold $(M,g)$ is called \textit{\hK} if it admits three complex structures $I,J$ and $K$ which:
\begin{enumerate}
\item are compatible with the metric, i.e. 
\[ g(\cdot,\cdot)=g(I\cdot,I\cdot)=g(J\cdot,J\cdot)=g(K\cdot,K\cdot) \]
\item  verify the quaternionic relations:
\[ IJ=-JI=K \]
\item are parallel with respect to the Levi-Civita connection given by $g$.
\end{enumerate} 
\end{de}

In particular, a \hK\ manifold is \Ka\ with respect to its fixed metric and any complex structure $aI+bJ+cK$, with $a, b$ and $c$ real constants verifying $a^2+b^2+c^2=1$.

Equivalently, we could say that a $4n$-dimensional Riemannian manifold $(M,g)$ is \hK\ iff its holonomy group is a subgroup of Sp$(n)$.

\begin{de}
A holomorphic 2-form on a complex manifold $M$,  $\omega\in H^0(M,\Omega_M^2)$,  is called a \textit{holomorphic symplectic form} if it is nondegenerate in the following sense:
\[v \lrcorner \omega_x=0 \Rightarrow v=0, \ \  \forall x\in M,  \forall v\in T^{1,0}_xM, \text{ where }\lrcorner \text{ is the contraction.}\]

We call a manifold admitting such a form a \textit{holomorphic symplectic manifold}.
\end{de}

In particular, a holomorphic symplectic manifold $(M,\omega)$ has even complex dimension $2m$ and $\omega^m$ is a nowhere vanishing holomorphic section of the canonical bundle $K_M$=det $\Omega_M^1$. Thus, $K_M$ is holomorphically trivial and $c_1(M)=0$.

It can be easily seen that, once we fix a complex structure on a \hK\ manifold $M$, say $I$, there exists a holomorphic symplectic form $\omega$ on $(M,I)$ defined by:
\[\omega(\cdot,\cdot)=g(J\cdot,\cdot)+ig(K\cdot,\cdot)  \]

Thus, a \hK\ manifold is a holomorphic symplectic manifold (but not in a canonical way). In the compact case, the converse is also true:

\begin{te}(\textbf{Beauville}, \cite{b})
Let $(M,I)$ be a compact \Ka\ manifold admitting a holomorphic symplectic form. Then, for any \Ka\ class $\al\in H^2(M,\RR)$, there exists a unique metric $g$ on M which is \Ka\ with respect to $I$, representing $\alpha$, so that $(M,g)$ is \hK. 

Moreover, the manifold $(M,I)$ admits a metric with holonomy exactly $\mathrm{Sp}(m$) iff it is simply connected and admits a unique holomorphic symplectic form up to multiplication by a scalar.  
\end{te}

\begin{ob}
The existence and uniqueness of the \Ka\ metric representing the given \Ka\ class comes from Yau's theorem: it is exactly the unique representative in the class that has vanishing Ricci curvature. Consequently, the holomorphic symplectic form in the theorem is parallel with respect to the Levi-Civita connection given by this Ricci-flat metric.
\end{ob}

\section{Twisted holomorphic symplectic manifolds}

We will now concentrate on the twisted case, and see that the situation is similar to the non-twisted one. Specifically, we will show that a \Ka\ manifold admitting a nondegenerate twisted holomorphic form admits a locally \hK\ metric which is moreover \Ka\ for the given complex structure. With respect to the connection induced by this metric, the form will be parallel.

\begin{de}
A Riemannian manifold $(M^{4m},g)$ is called \textit{locally \hK } if its universal cover with the pullback metric is \hK\ or, equivalently, if the restricted holonomy group $\text{Hol}^0(g)$ is a subgroup of Sp$(m$). %
\end{de}

\begin{de}
We will call a \Ka\ manifold $(M,I,g)$ \textit{\Ka\ locally \hK\ (Klh)} if it admits a rank 2 real sub-bundle of $\text{End}_\RR(TM)$ that has, around each point, a local parallel frame consisting in two orthogonal complex structures which, together with $I$, verify the quaternionic relations. %
\end{de} 

\begin{ob} Notice that a Riemannian manifold $(M,g)$ is locally \hK\ exactly when there exists a rank 3 real sub-bundle $\mathcal{G}$ of $\text{End}_\RR(TM)$ which admits, around each point, a local parallel frame consisting in three orthogonal complex structures verifying the quaternionic relations. When, moreover, $\mathcal{G}$ admits as a global section a complex structure $I$ with respect to which $g$ is \Ka, the manifold is Klh.
\end{ob}

\begin{de}
A compact manifold $(M,I,L,\omega)$ is called \textit{twisted holomorphic symplectic} if $I$ is a complex structure on $M$ and there exists  a holomorphic line bundle $L$ over $M$ and a nondegenerate $L$-valued holomorphic form:
$$\omega\in H^0(M,\Omega^2_M\otimes L).$$
 
\end{de}

\begin{ob}
The existence of a twisted holomorphic symplectic form implies that $M$ is of even complex dimension $2m$. Moreover, $\omega^m$ is a nowhere vanishing holomorphic section of the line bundle $K_M\otimes L^m$. Thus, we have that $L^m\cong K_M^*$. In particular, any metric on $M$ naturally induces one on $L$, and we also have
$$c_1(M)=mc_1(L).$$ 

\end{ob}

\begin{ob} 
Any compact complex surface $M$ is twisted-symplectic in a tautological way. Simply take $L$ to be $K_M^*$, so that
$$H^0(M,\Omega^{2}_M\otimes L)=H^0(M,K_M\otimes K_M^*)=\CC.$$
Thus any non-zero complex number represents a twisted-symplectic form. Therefore, the class of twisted holomorphic symplectic manifolds is interesting only starting from complex dimension $4$. 
\end{ob}

Our main result in this section is the following:

\begin{te}
Let $(M^{2m},I,L,\omega)$, $m>1$, be a compact twisted holomorphic symplectic manifold of \Ka\ type, and let $\al\in H^2(M,\RR)$ be a \Ka\ class. Then there exists a unique \Ka\ metric $g$ with respect to $I$ representing $\al$ so that $(M,g,I)$ is Klh. Moreover, $L$ is flat and $\omega$ is parallel with respect to the natural connection induced by $g$ on $L$.
\end{te}

\begin{demo}
Let $\{U_i\}$ be a trivializing open cover for the line bundle $L$ and $\sigma_i \in H^0(U_i,L)$ be holomorphic frames, so that the holomorphic transition functions $\{g_{ij}\}$ are given by  $\sigma_i = g_{ij}\sigma_j$. Then, if we write over $U_i$
$$\omega=\omega_i\otimes\sigma_i$$
we get local holomorphic symplectic forms $\omega_i$ that verify, on $U_i \cap U_j$, 
$\omega_i = g_{ji}\omega_j$.

The $\omega_i$'s, being holomorphic, induce the morphisms of sheaves of $\Oo_{U_i}$-modules over $U_i$:
$$L_i:\Omega^k_{U_i}\rightarrow\Omega^{k+2}_{U_i}$$
$$L_i=\omega_i\wedge\cdot$$

\begin{lm} 
For $m>1$, shrinking the $U_i$'s if necessary, we have an isomorphism of sheaves of $\Oo_{U_i}$-modules:
$$\Omega^3_{U_i}\cong\Omega^1_{U_i}\oplus\Omega^3_{0,U_i}$$
where $\Omega^3_{0,U_i}$ is the sheaf $\ \ker (L_i^{m-2}:\Omega^3_{U_i}\rightarrow\Omega^{n-1}_{U_i})$ and $n=2m$.

\end{lm}

\begin{demo}We claim that $L_i^{m-1}:\Omega^1_{U_i}\rightarrow\Omega^{n-1}_{U_i}$ is an isomorphism of sheaves over $U_i$. We inspect this at the germ level, so we fix $z\in U_i$. Since the corresponding free $\Oo_z$-modules have the same dimension, it suffices to prove the injectivity of $L^{m-1}_{i,z}$. But this becomes a trivial linear algebra problem, noting that we can always find a basis over $\CC$ in $T^{1,0}M_z^*$ $\{e_1,\ldots,e_m,f_1,\ldots, f_m\}$ so that
$$\omega_i(z)=\sum_{s=1}^me_s\wedge f_s.$$

Next, from
$$\xymatrix{
&\Omega^3_{U_i} \ar[rd]^{L_i^{m-2}}\\
&\Omega^1_{U_i}\ar[u]^{L_i} \ar[r]^{L_i^{m-1}}  &\Omega^{n-1}_{U_i}} $$
we get that $L_i$ is injective and $L_i^{m-2}$ is surjective. Hence, we have an exact sequence of sheaves: 
$$
\xymatrix{
0 \ar[r] &\Omega^3_{0,U_i} \ar[r] &\Omega^3_{U_i} \ar[r]^{L_i^{m-2}} &\Omega^{n-1}_{U_i}\ar[r] &0}
$$
Now, we can take the $U_i$'s small enough so that all the sheaves above are free. Thus, the sequence splits and identifying $\Omega^{n-1}_{U_i}\cong\Omega^1_{U_i}$ with a subsheaf of $\Omega^3_{U_i}$ via $L_i$,  we get the desired isomorphism. This ends the proof of the lemma.\end{demo}

Now, we have $d\omega_i\in\Omega^3_M(U_i)$, so, refining the cover $\{U_i\}$ if necessary,  we can write:
\begin{equation}
d\omega_i=\omega_i\wedge\theta_i+\xi_i
\end{equation}
with $\theta_i\in\Omega^1_M(U_i)$ and $\xi_i\in\Omega_{0,M}^3(U_i)$ holomorphic sections uniquely determined by the previous lemma. Since $\omega_i=g_{ji}\omega_j$, we get: 
$$dg_{ji}\wedge\omega_j+g_{ji}d \omega_j=g_{ji} \omega_j \wedge \theta_i+\xi_i$$
whence
$$\omega_j \wedge \theta_j+ \xi_j=d\omega_j = \omega_j\wedge\theta_i+\frac {1}{g_{ji}}\xi_i-\frac{dg_{ji}}{g_{ji}}\wedge\omega_j$$

Thus the $\theta_i$'s change by the rule:
\begin{equation}
\theta_i=\theta_j+d\log g_{ji}.
\end{equation}

Now, on a compact \Ka\ manifold we have $H^1(M,\Omega_M^1)\cong H^{1,1}(M,\CC)\subset H^2(M,\CC)$ and $c_1(L)\in H^{1,1}(M,\CC)$ is represented, via this isomorphism, by the \v{C}ech cocycle $\{d\log g_{ji}\}_{ji}\in \check{Z}^1(\{U_i\}_i, \Omega_M^1)$. On the other hand, from (2.2) we have that $\delta(\{-\theta_i\}_i)=-(\theta_j-\theta_i)=d\log g_{ji}$, where $\{-\theta_i\}_i\in \check{Z}^0(\{U_i\}_i,\Omega_M^1)$ and $\delta$ is the differential in the \v{C}ech complex. Therefore, $\{d\log g_{ji}\}_{ji}$ is exact, hence $c_1(L)=0$.

Thus we also get $c_1(M)=mc_1(L)=0$. So, by Yau's theorem, there exists a unique Ricci-flat \Ka\ metric $g$ whose fundamental form $\omega_g$ represents the given class $\al$. 

Now, on $\Omega^{2,0}_M\otimes L$ we have the Weitzenb\"{o}ck formula (see for instance \cite{m}):
$$2\cp^*\cp =\nabla^*\nabla+\mathcal{R}$$
where $\nabla$ is the naturally induced connection by $g$ on $\Omega^{2,0}_M\otimes L$ and  $\mathcal{R}$ is a curvature operator which on decomposable sections is given by:
$$\mathcal{R}(\beta\otimes s)=i\rho_g\beta\otimes s + \beta\otimes \text{Tr}_{\omega_g}(i\Theta(L))s$$
with $\rho_g:\Omega^{2,0}_M\rightarrow \Omega^{2,0}_M$ the induced action of the Ricci form on $\Omega^{2,0}_M$. Now, since $g$ is Ricci-flat, $\rho_g\equiv 0$. Also, if we consider the curvatures induced by $g$, we have:
$$0=-i\rho=\Theta(K^*_M)=\Theta(L^m)$$
so the induced connection on $L$ is flat and $\mathcal{R}$ vanishes.

Hence, applying the Weitzenb\"{o}ck formula to $\omega$, we get $0=\nabla^*\nabla\omega$ or also, after integrating over $M$, $\|\nabla\omega\|^2_{L^2}=0$. Thus $\nabla\omega=0$.

Finaly, if we let $\pi:(\tilde{M},\tilde{g},\tilde{I})\rightarrow(M,g,I)$ be the universal cover with the pullback metric and complex structure, we have that $\pi^*L$ is trivial and $\tilde\omega=\pi^*\omega\in H^0(\tilde{M},\Omega_{\tilde M}^2)$ is a holomorphic symplectic form. By the Cheeger-Gromoll theorem, $\tilde M\cong \CC^l\times M_0$, where $M_0$ is compact, simply connected, \Ka, Ricci-flat, and $\CC^l$ has the standard \Ka\ metric. Moreover, by the theorems of de Rham and Berger, the holonomy of $M_0$ is a product of groups of type Sp$(k)$ and SU$(k)$.  We have that $\tilde\omega$ is a parallel section of 
$$\textstyle\bigwedge^2T^*\tilde M= \textstyle\bigwedge^2 \text{pr}_1^*T^*\CC^l\oplus (\text{pr}_1^*T^*\CC^l\otimes \text{pr}_2^* T^*M_0)\oplus \textstyle\bigwedge^2 \text{pr}_2^*T^*M_0$$ 
But $\text{pr}_1^*T^*\CC^l\otimes \text{pr}_2^* T^*M_0\cong (T^*M_0)^{\oplus l}$ has no parallel sections by the holonomy principle, so $\tilde\omega$ is of the form 
$\omega_c+\omega_0$, with $\omega_c$, $\omega_0$ holomorphic symplectic forms on $\CC^l$, $M_0$ respectively. Thus, $l$ is even, so $\CC^l$ is \hK, and also, by Theorem 1.3, $M_0$ is \hK. It follows that  $(M,g,I)$ is Klh.

This concludes the proof of the theorem.\end{demo}

\begin{ob}
Notice that the relation (2.2) says that the $\theta_i$'s behave like connection forms on $L$ over $U_i$, i.e. the differential operator $D:C^\infty(M,L)  \rightarrow C^\infty(M, T^*M\otimes L)$ given over $U_i$ by $D(f\otimes \sigma_i)=(df-\theta_i)\otimes \sigma_i$ is a well defined connection on $L$.
 We will then have that its curvature $\Theta(D)_{U_i}=-d\theta_i$ is of type (2,0) and represents the first Chern class of $L$ up to a constant. But it is well known that on a compact \Ka\ manifold, the image of $c_1$ in $H^2(M,\CC)$ consists of the integral classes from $H^{1,1}(M,\CC)$, so $c_1(L)=0$. Note that $D^{1,0}$ is actually a holomorphic connection on $L$.
\end{ob}

\begin{ob}
The flat connection induced by $g$ on $L$ does not depend on the \Ka\ class $\al$. It is uniquely determined by $\omega$ and is equal to the connection given in the previous remark. To see this, let $D^g$ be the Chern connection on $L$ induced by $g$ and write $D^g\sigma_i=\tau_i\otimes\sigma_i$. Then we have:
$$0=\nabla\omega=\nabla\omega_i\otimes\sigma_i +\omega_i\otimes\tau_i \otimes\sigma_i.$$

So, denoting by $a:
\Omega_M^2\otimes T^*M^c\otimes L \rightarrow(\Omega_M^{3,0}\oplus\Omega_M^{2,1})\otimes L$ the antisymmetrization map, we get:
$$d\omega_i=a(\nabla\omega_i)= -\omega_i\wedge\tau_i.$$
Thus, by (1.1) we deduce that $\xi_i=0$ and $\tau_i=-\theta_i$, i.e. $D^g=D$.

\end{ob}

\begin{ob}
If we only suppose that $\omega$ is a non degenerate $(2,0)$ twisted form, not necessarily holomorphic, then $\omega$ still induces a connection on $L$ in the same manner. This time, $L_i:\Omega^{k,0}_{U_i}\rightarrow\Omega^{k+2,0}_{U_i}$ are morphisms of sheaves of $\E_{U_i}$-modules, and induce isomorphisms $\Omega^{3,0}_M(U_i)\cong\Omega^{1,0}_M(U_i)\oplus\Omega^{3,0}_{0,M}(U_i)$.  Writing  $$\Omega^{3,0}_M(U_i)\ni\partial \omega_i=\omega_i\wedge\theta_i+\xi_i,$$ we get the $(1,0)$ forms $\theta_i$ which define a connection $D$ just as in Remark 2.9. It is only at this point that the holomorphicity of $\omega$ becomes essential in order to have that $D$ defines a holomorphic connection on $L$. 
\end{ob}

\section{A characterization}

In this section, we want to investigate the converse problem. It is not true that all Klh manifolds are twisted holomorphic symplectic. Already we will see that a product of strictly twisted symplectic manifolds is never twisted symplectic, but it turns out that being reducible is not the only obstruction. In what follows, we will give some description of twisted holomorphic symplectic manifolds and their fundamental groups. 

By a strictly twisted holomorphic symplectic manifold we always mean a twisted holomorphic symplectic manifold $(M,I,L,\omega)$ such that the line bundle $L$ is not holomorphically trivial.

\begin{pr}
A compact \Ka\ manifold $M$ of complex dimension $>2$ is twisted holomorphic symplectic iff there exists a holomorphic symplectic form $\omega_0$ on its universal cover $\tilde M$ so that the action of  $\  \Gamma=\pi_1(M)$ on $H^0(\tilde M,\Omega_{\tilde M}^2)$ preserves $\CC\omega_0$. In particular, any twisted-symplectic manifold is a finite cyclic quotient of a \hK\ manifold.
\end{pr}

\begin{demo}
Suppose first that $M$ admits a twisted-symplectic form
$$\omega\in H^0(M,\Omega^2_M\otimes L).$$
Then, by Theorem 2.7,  $L$ is flat, and thus given by a unitary representation $\rho:\Gamma\rightarrow$U$(1)$, i.e. if we see $\pi:\tilde M \rightarrow M$ as a $\Gamma$-principal bundle over $M$, we have $L=\tilde M \times_\rho \CC$.

Let $s_i:U_i\rightarrow\tilde M$ be local sections of  $\pi:\tilde M \rightarrow M$ over a trivializing cover $\{U_i\}$. We then have  $s_i=\gamma_{ij}s_j$ on $U_i \cap U_j$, where $\gamma_{ij}: U_i\cap U_j \rightarrow \Gamma$ are the transition functions for $\tilde M$. Then, $\sigma_i:=[s_i,1]$ are local frames for $L$, where $[\cdot,\cdot]$ denotes the orbit of an element of $\tilde M\times \CC$ under the left action of $\Gamma$. The locally constant functions $g_{ij}:=\rho(\gamma_{ij}^{-1})$  are the transition functions for $L$ verifying 
$$\sigma_i=[\gamma_{ij}s_j,1]=[s_j,\rho(\gamma_{ij}^{-1})]=g_{ij}\sigma_j.$$

Since $\pi^*L$ is trivial, there exist $f_i\in\Oo^*_{\tilde M}(\pi^{-1}U_i)$ such that $\pi^*g_{ij}=\frac{f_i}{f_j}$ on $\pi^{-1}U_i\cap\pi^{-1}U_j$. Also, the sections $\frac{\pi^*\sigma_i}{f_i}\in H^0(\pi^{-1}U_i,\pi^*L)$ all coincide on intersections and are non vanishing, thus giving a global frame for $\pi^*L$ which we can suppose equal to 1, so that $\pi^*\sigma_i=f_i$. Thus, if we write $\omega=\omega_i\otimes\sigma_i$ and define $\omega_0:=\pi^*\omega$, we get:
$$\omega_0|_{\pi^{-1}U_i}=\pi^*\omega_if_i$$
and, for any $\gamma\in\Gamma$:
$$\gamma^*\omega_0|_{\pi^{-1}U_i}=\pi^*\omega_i\gamma^*f_i=\omega_0\frac{\gamma^*f_i}{f_i}$$
Moreover, for any $\gamma$, we have on $\pi^{-1}U_i\cap\pi^{-1}U_j$ :
$$\frac{f_j}{f_i}=\frac{f_j\circ\gamma}{f_i\circ\gamma} \Leftrightarrow g_{ij}\circ\pi=g_{ij}\circ\pi\circ\gamma$$
hence the constant function $\frac{f_i\circ \gamma}{f_i}$ does not depend on $i$. 

On the other hand, we have:
\begin{equation}
\frac{\gamma^*f_i}{f_i}=\frac{[s_i\circ\pi\circ\gamma,1]}{[s_i\circ\pi,1]}=\frac{[s_i\circ\pi,\rho(\gamma^{-1})]}{[s_i\circ\pi,1]}=\frac{1}{\rho(\gamma)}
\end{equation}

Hence $\Gamma$ preserves the subspace $\CC\omega_0\subset H^0(\tilde M,\Omega_{\tilde M}^2)$ and  $\rho$ is determined by the action of $\Gamma$ on the holomorphic symplectic form $\omega_0$ by:
 $$\frac {1} {\rho(\gamma)}\cdot \omega_0=\gamma^*\omega_0.$$

 Conversely, suppose a holomorphic symplectic form $\omega_0$ is an eigenvector for $\Gamma$ acting on $H^0(\tilde M,\Omega_{\tilde M}^2)$. Define 
 $$\rho:\Gamma\rightarrow \CC^*$$
$$\gamma\mapsto \frac{\omega_0}{\gamma^*\omega_0}$$
Let $L:=\tilde{M}\times_{\rho}\CC$ and, with the same data for $L$ as before, define $\omega\in H^0(M,\Omega_M^2\otimes L)$  by $\omega|_{U_i}=\omega_i\otimes\sigma_i$, where $\omega_i=s_i^*\frac{\omega_0}{f_i}$. Then $\omega$ is twisted holomorphic symplectic and, seeing $s_i\pi$ as an element of $\Gamma$, we have, by (3.1), on $\pi^{-1}(U_i)$ : 
$$\pi^*\omega= \frac{\pi^*s_i^*\omega_0}{\pi^*s_i^*f_i}f_i=\frac{1}{\rho(s_i\pi)}\omega_0\frac{f_i}{(s_i\pi)^*f_i}=\omega_0.$$

To prove the last part, suppose $M$ is twisted-symplectic and let, as in Theorem 2.7,  $\tilde M=\CC^{2l}\times M_1\times\ldots M_k$, with $M_i$ irreducible \hK\ manifolds. The manifold $M$ has a finite \'etale cover $M'=\TT^{2l}\times M_1\times\ldots M_k$ so that  $M=M'/ \Gamma'$ and  $\Gamma\cong \ZZ^{4l}\ltimes \Gamma'$. The symplectic form $\omega_0$ is preserved under the action of $\ZZ^{4l}$, so it descends to a holomorphic symplectic form on $M'$, which we will also denote by $\omega_0$. The group $\Gamma'$ preserves $\CC\omega_0$.

Let $\rho':\Gamma'\rightarrow$U(1) be the representation induced by $\rho$. Denote by $N'$ its kernel, and by $N:=\ZZ^{4l}\ltimes N'$. Then $N$ is normal inside $\Gamma$, so there exists a Galois covering $M_N\rightarrow M$ with $\pi_1(M_N)=N$. Moreover, since $\pi_1(M')=\ZZ^{4l}$ is normal in $N$, also $M'\rightarrow M_N$ is a covering whose deck transformation group is $N/\ZZ^{4l}\cong N'$. 

We thus have that $M_N\cong M'/N'$ and $N'$ preserves $\omega_0$, so $\omega_0$ descends to $M_N$. Since $M_N$ is compact holomorphic symplectic, it is \hK.

Finally, $\rho(\Gamma)=\rho'(\Gamma')$ is a finite subgroup of U(1), so cyclic, and $\Gamma/N\cong\Gamma'/N'\cong\rho(\Gamma)$, so $M_N$ is a finite cyclic covering of $M$.

This concludes the proof of the proposition.
\end{demo}

\begin{co}
A compact strictly twisted holomorphic symplectic manifold of dimension $>2$ is de Rham irreducible.
\end{co}

\begin{demo}
Suppose $M\cong M_1\times M_2$ is strictly twisted holomorphic symplectic. Let $\tilde{M}\cong\tilde{M_1}\times\tilde{M_2}$ be a finite \'etale cover of $M$ with holomorphic symplectic form $\omega_0=\omega_1+\omega_2$ preserved up to constants by $\Gamma'\cong\Gamma_1'\times \Gamma_2'$, where $\pi_1(M_i)=\ZZ^{2l_i}\ltimes\Gamma_i'$, $i=1,2$. 

Then we should have that $\rho(\Gamma')=\rho(\Gamma_1')\times\rho(\Gamma_2')$ is a non-trivial cyclic group of the same order as $\rho(\Gamma_1')$, $\rho(\Gamma_2')$, which is impossible.
\end{demo}

\begin{co}
A compact locally irreducible \Ka\ manifold of dimension $>2$ is Klh iff it is twisted holomorphic symplectic. In this case, the twisted-symplectic form is valued in the canonical bundle. 
\end{co}

\begin{demo}
Let $M$ be a locally irreducible Klh manifold and $\tilde M$ its universal cover endowed with a holomorphic symplectic form $\omega_0$. Since $\tilde M$ is irreducible, it is compact and $H^0(\tilde M,\Omega_{\tilde M}^2)=\CC\omega_0$. Hence $\Gamma=\pi_1(M)$ preserves $\CC\omega_0$ in a trivial way and $M$ is twisted holomorphic symplectic by the previous proposition. 

In particular, this implies that $\Gamma$ is cyclic. Let $d$ be its order. Then $d|m+1$, where dim$M=2m$. To see this, let $\gamma\in\Gamma$ be a generator, so that $\gamma^*\omega_0=\xi\cdot\omega_0$, with $\xi$ a primitive $d$-root of unity. Since $\gamma$ has no fixed points, by the holomorphic Lefschetz fixed-point formula we must have that its Lefschetz number, which by definition is:
$$L(\gamma)=\sum_q (-1)^q\tr\gamma^*|_{H^q(\tilde M,\Oo)}$$
must vanish. On the other hand, we have 

$$\widebar{H}^*(\tilde M,\Oo_{\tilde M})\cong H^0(\tilde M,\Omega^*_{\tilde M})\cong\frac{\CC[\omega_0]}{(\omega_0^{m+1})}$$
so $L(\gamma)=1+\xi+\ldots+\xi^m$. Thus, $L(\gamma)=0$ implies $d|m+1$.

Let $\rho:\Gamma\rightarrow$U(1) be given by the action of $\Gamma$ on $\omega_0$ and $L:=\tilde M\times_{\widebar{\rho}}\CC$, so that the twisted holomorphic symplectic form is $L$-valued. Since the action of $\Gamma$ on $K_{\tilde M}$ is given by $\rho^m$, we also have that $K_M=\tilde M\times_{\rho^m}\CC$. Now, $\rho^{m+1}=1$ implies $\widebar{\rho}^m\cdot\widebar{\rho}=1$, or also $K_M^*\otimes L=\underline{\CC}$, i.e. $L\cong K_M$. \end{demo}

\begin{ob}
For a twisted holomorphic symplectic manifold $(M,I,L,\omega)$, we always have, by Remark 2.5, that $L$ is a root of $K_M^*$. In the particular case when $M$ is locally irreducible, we obtain, moreover, that $L$ (and thus also $K_M$) is a torsion element of the Picard group, and that $L$ is precisely (up to isomorphism) $K_M$.
\end{ob}

It is difficult to give a nice criterion for being twisted holomorphic symplectic in the case of de Rham irreducible, locally reducible \hK\ manifolds. We can, though, give a somewhat more precise description of fundamental groups of twisted holomorphic symplectic manifolds. For this, we first give some lemmas concerning isometries of Riemannian products.

\begin{lm} Let $(M_i,g_i)$ be complete locally irreducible Riemannian manifolds of dimension >1 and let $M_0=M_1\times\ldots\times M_k$ be endowed with the product metric.  Let $\gamma$ be an isometry of $M_0$ and let $\gamma_i:=p_i\gamma$, where $p_i:M_0\rightarrow M_i$ are the canonical projections. Then $\gamma_i$ is of the form $\gamma_i=\tilde \gamma_ip_{\sigma(i)}$, where $\tilde\gamma_i:M_{\sigma(i)}\rightarrow M_i$ is an isometry and $\sigma$ a permutation of $\{1,\ldots,k\}$.
\end{lm}

\begin{demo} We have that $\tilde g_i:=\gamma_i^* g_i$ is a parallel section of $S^2(T^*M_0)$. On the other hand, 
$$S^2(T^*M_0)\cong \sum S^2(T^*M_i)\oplus \sum_{i<j}(T^*M_i\otimes T^*M_j).$$

Now, $T^*M_i\otimes T^*M_j$ admits no parallel section for $i<j$, while the space of parallel sections of $S^2(T^*M_i)$ is exactly $\RR g_i$. Indeed, by the holonomy principle, this is equivalent to saying that $G_i\times G_j$ has no fixed points when acting on $T_x^*M_i\otimes T_y^*M_j$, while the only  $G_i$-invariant elements of $S^2(T_x^*M_i)$ are the multiples of $(g_i)_x$, where $x\in M_i$, $y\in M_j$ are any points and $G_s$ is the restricted holonomy group of $M_s$, $s=1,...,k$.  The first assertion follows from the dimension hypothesis and the more general fact that if $U$ is a $G$-irreducible space and $V$ a $H$-irreducible space, then $U\otimes V$ is a $G\times H$-irreducible space. The second assertion is equivalent to Schur's lemma if we identify $S^2(T_x^*M_i)$ with the symmetric endomorphisms of $T_x^*M_i$ via $g_i$.

Next, we want to show that for every $i$, there is exactly one $j=j(i)$ so that $a_{ij}\neq 0$. Thus, if we let $A(i)=\{j|a_{ij}\neq 0\}$, we need to show that $A(i)\neq\emptyset$ for each $i$ and $A(i)\cap A(j)=\emptyset$ for all $i\neq j$. Now, since $g_i$ is definite and $d\gamma_i$ is surjective, we have that $\ker \tilde g_i:=\{X \in TM_0| \tilde g_i(X,\cdot)=0\}=\ker d\gamma_i$. Hence, since $\ker d\gamma_i\neq TM$, the first assertion follows.

For the second assertion, first note that $\ker\tilde g_i\cap TM_k\neq 0$ iff $a_{ik}=0$, in which case $TM_k\subset \ker \tilde g_i$. Therefore, $(\ker d\gamma_i)^\perp=\sum_{j\in A(i)} TM_j$. %
Hence, for $i\neq j$, $A(i)\cap A(j)=\emptyset$ is equivalent to $\{0\}=(\ker d\gamma_i)^\perp\cap(\ker d\gamma_j)^\perp=(\ker d\gamma_i + \ker d\gamma_j)^\perp$. But we have
$$\ker d\gamma_i + \ker d\gamma_j=d\gamma^{-1}(\ker d p_i+ \ker dp_j)=d\gamma^{-1}(\sum_{s\neq i} TM_s+\sum_{s\neq j}TM_s)=TM.$$
 
It follows that there exists a permutation $\sigma$ of $\{1,\ldots,k\}$ so that $A(i)=\{\sigma(i)\}$ for each$i$. Hence, since for any $j$, $\sum_i a_{ij}=1$, we have that $a_{i\sigma(i)}=1$ and $\gamma_i=\tilde\gamma_ip_{\sigma(i)}$ with $\tilde\gamma_i:M_{\sigma(i)}\rightarrow M_i$ an isometry.
\end{demo}

In what follows, we will omit writing the projections and identify $\gamma_i$ with $\tilde\gamma_i$.

\begin{lm}
Let $M_i$ be irreducible compact \hK\ manifolds, and $M_0=M_1\times\ldots\times M_k$ be endowed with the product metric and a holomorphic symplectic form $\omega_0$. Then any isometry of $M_0$ preserving $\omega_0$ has fixed points.
\end{lm} 

\begin{demo}
By Theorem 1.3, the $M_i$'s are simply connected and admit unique holomorphic symplectic forms $\omega_i$ up to a scalar, so we have:
$$H^0(M_0,\Omega_{M_0}^2)=\CC\omega_1\oplus\ldots\oplus\CC\omega_k.$$
Hence we can suppose, after rescaling the $\omega_i$'s, that $\omega_0=\omega_1+\ldots+\omega_k$. Let $\gamma$ be an isometry of $M_0$ with $\gamma^*\omega_0=\omega_0$.

Consider first the case where all $M_i$ are isometric, so that $M_0\cong M_1^k$.  Let $\sigma$ be the permutation determined by $\gamma$ as in the previous lemma and let $l$ be the order of $\sigma$. If we define, for $i=1,...,k$: 
$$\gamma'_i=\gamma_i\gamma_{\sigma(i)}\ldots\gamma_{\sigma^{l-1}(i)}$$
then $\gamma^l(x_1,\ldots,x_k)=(\gamma'_1(x_1),\ldots,\gamma'_k(x_k))$. If $\gamma$ acts freely, then also $\gamma^l$ acts freely. Otherwise, suppose $\gamma^l(y_1,\ldots,y_k)=(y_1,\ldots,y_k)$, let $i_1,\dots ,i_t \in \{1,\ldots ,k \}$ represent the orbits of $<\sigma>$ and define $(x_1,\ldots,x_k)$ by
$$x_{i_j}:=y_{i_j}, \ x_{\sigma(i_j)}:=\gamma_{\sigma(i_j)}\ldots\gamma_{\sigma^{l-1}(i_j)}(y_{i_j}), \ \ldots, \  x_{\sigma^{l-1}(i_j)}:=\gamma_{\sigma^{l-1}(i_j)}(y_{i_j})$$
The fact that $\gamma_{i_j}\gamma_{\sigma(i_j)}\ldots\gamma_{\sigma^{l-1}(i_j)}(y_{i_j})=y_{i_j}$ implies that $(x_1,\ldots,x_k)$ is a fixed point for $\gamma$, contradiction.

Now, $\gamma^*\omega_0=\omega_0$ implies $(\gamma^l)^*\omega_0=\sum_i(\gamma_i')^*\omega_1=\omega_0$, or also $(\gamma'_i)^*\omega_1=\omega_1$ for any $i=1,...,k$. On the other hand, the fact that $\gamma^l$ acts freely implies that some $\gamma'_{i_0}$ acts freely on $M_1$. By the holomorphic Lefschetz fixed-point formula, its Lefschetz number must then vanish. But $L(\gamma'_{i_0})=m+1$, where dim$M_1=2m$, contradiction.

In the general situation, write $M_0=(M_1)^{k_1}\times\ldots\times (M_s)^{k_s}$, with $M_i$ irreducible and $M_i\ncong M_j$  for all $i\neq j$. By the previous lemma, $\gamma=(\gamma_1,\ldots, \gamma_s)$, with $\gamma_i$ an isometry of $(M_i)^{k_i}$. Again, $\gamma^*\omega_0=\omega_0$ implies $\gamma_i^*\omega_{i0}=\omega_{i0}$, where the $\omega_{i0}$'s are the induced symplectic forms on $(M_i)^{k_i}$, $i=1,\ldots,s$. Also, if $\gamma$ acts freely on $M_0$, then some $\gamma_i$ acts freely on $(M_i)^{k_i}$ and we already showed that this is impossible.
\end{demo}

\begin{ob}
We can now say slightly more about twisted holomorphic symplectic manifolds $M$ with compact universal cover $\tilde M$. In this case, with the notations of Proposition 3.1, $\Gamma=\Gamma'$, the representation $\rho$ is faithful by the previous lemma, so $\Gamma=\rho(\Gamma)$ is cyclic. Thus, if $\gamma$ is a generator of $\Gamma$ of order $d$ and $\gamma^*\omega_0=\xi\omega_0$, then $\xi$ is necessarily a primitive $d$-root of unity. Moreover, if we write $\gamma=(\gamma_1,\ldots ,\gamma_k)$ just as in Lemma 3.5, then all $\gamma_i$'s must have the same order $d$. To see this, let $d_i=$ord$\gamma_i$. Then $d_i|d=$lcm$(d_i)_i$. Since $\gamma^*\omega_0=\sum_i\gamma_i^*\omega_i=\xi\sum_ip^*_i\omega_i$, we have, for all $i$, $\gamma_i^*\omega_i=\xi\omega_{\sigma(i)}$, hence $\xi^{d_i}=1$. But $\xi$ was primitive, so $d_i=d$. We can conclude:
\end{ob}

\begin{co}
If the fundamental group of a compact twisted holomorphic symplectic manifold is finite, then it is cyclic and of the form $\Gamma=<\gamma=(\gamma_1,\ldots,\gamma_k)>$, with $\gamma_i$ isometries of the irreducible components of the universal cover, all of the same order.  \end{co}

\begin{ob}
When $M$ is twisted holomorphic symplectic but $\tilde M$ is not compact, it is not necessarily the case for $\Gamma'$ to be cyclic, i.e. $\rho':\Gamma'\rightarrow$U(1) need not be faithful. By the same type of arguments as in Lemma 3.5, it can be seen that an element of $\Gamma'$ is of the form $\gamma=(\gamma_T,\gamma_0)$, with $\gamma_T\in$Aut$(\TT^{2l})$ and $\gamma_0\in$Aut$(M_0)$.  There exist fixed point free complex symplectomorphisms of $\TT^{2l}$ of finite order (for instance translation by a torsion element $a \in\TT^{2l}$ ). So, if $\gamma_T$ is one and $\gamma_0$ is a symplectomorphism of $M_0$ of the same order as $\gamma_T$, $(\gamma_T,\gamma_0)$ is an element in the kernel of $\rho'$.

\end{ob}

\begin{co}
A compact strictly twisted holomorphic symplectic manifold $M$ of dimension $>2$ with finite fundamental group is projective.
\end{co}

\begin{demo}
Let $\pi:\tilde M\rightarrow M$ be the compact universal covering, where, by Theorem 2.7, $\tilde M=M_1\times\ldots\times M_k$ with $M_i$ irreducible \hK\ manifolds. Then, by Lemma 3.6, each $M_i$ admits an automorphism which is not symplectic. By a result of \cite{bs}, such manifolds are necessarily projective, hence so is $\tilde M$. But it is a well known fact that a compact \Ka\ manifold is projective if and only if some finite unramified covering is, thus the conclusion follows. \end{demo}

\section{Final remarks}

\begin{ob}
Concerning examples, finding locally irreducible Klh manifolds is equivalent to finding a fixed point free automorphism $\gamma$ of an irreducible symplectic manifold, so that all powers of $\gamma$ also act freely. 

In complex dimension 2, by Remark 2.6, all manifolds are twisted holomorphic symplectic. On the other hand, all compact Klh surfaces are either tori or Enriques surfaces, which are quotients $K/<\iota>$, with $K$ a K3 surface admitting a fixed point free involution $\iota$. 

Next, out of a K3 surface $(K, \iota)$ as before, one can construct a twisted holomorphic symplectic manifold of any even dimension $2m$. Simply take $K^{m}/<\iota,\ldots,\iota>$. 

To find twisted holomorphic symplectic manifolds with higher order of the fundamental group, one needs to look at other irreducible \hK\ manifolds. For the Hilbert schemes of points on K3 surfaces,  see \cite{b} for the construction, all known automorphisms have fixed points, so we have no hope of constructing examples out of them.  

On the other hand, there is hope with the generalized Kummer varieties $K_r$, see again \cite{b} for the definition. In \cite{bns} and \cite{o} the authors find fixed point free cyclic groups of automorphisms $\Gamma$ of order 3 for the manifolds $K_2$ and $K_5$, and of order 4 for $K_3$. Taking the corresponding quotients give the desired examples of locally irreducible twisted holomorphic symplectic manifolds of dimension 4, 10 and 6, respectively.
\end{ob}

\begin{ob}
In order to actually classify twisted holomorphic symplectic manifolds, one should be able to classify fixed point free groups of automorphisms of irreducible \hK\ manifolds. The problem is clear in low dimension. It is also clear that if the Hilbert schemes of points on K3 surfaces admit such groups, then the corresponding automorphisms are not natural, i.e. do not arise from automorphisms of the K3 surface. On the other hand, for the moment there are no known exemples of such non-natural automorphisms. For the generalized Kummer varieties, there exist some examples formed out of natural automorphisms, but we do not know a classification of such groups.
\end{ob}

\begin{ob}
The \Ka\ hypothesis was heavily used to show that twisted holomorphic symplectic manifolds are locally \hK, particularly when applying Yau's theorem. On the other hand, we have no examples of twisted holomorphic symplectic manifolds in the non-\Ka\ case, and the problem is open in this generality. 
\end{ob}

\begin{ob}
Another direction to go would be to study the problem in the non holomorphic case, that is to be able to say which compact \Ka\ manifolds admit a non degenerate (2,0)-form valued in a complex line bundle. These are the manifolds admitting a topological Sp($m$)U(1) structure. As was mentioned in Remark 2.11, Theorem 2.7 does not hold without the holomorphicity assumption. A counterexample is given by the quadric $\QQ_6=$SO$(7)/$U$(3)\subset\PP^7\CC$, which is a \Ka\ manifold with topological Sp($m$)U(1) structure, see \cite{m2}, but is not Klh, since it has positive first Chern class.  
\end{ob}

\bigskip
\bigskip
\begingroup
\let\clearpage\relax

\endgroup

\bigskip

{\small \textsc{Nicolina Istrati, Univ Paris Diderot, Sorbonne Paris Cité, Institut de Mathématiques de Jussieu-Paris Rive Gauche, UMR 7586, CNRS, Sorbonne Universités, UPMC Univ Paris 06, F-75013, Paris, France} 

\textit{E-mail address}: {\scriptsize nicolina.istrati@imj-prg.fr}}

\end{document}